\documentclass{birkmult}
 \newtheorem{thm}{Theorem}[section]
 \newtheorem{cor}[thm]{Corollary}
 \newtheorem{lem}[thm]{Lemma}
 \newtheorem{prop}[thm]{Proposition}
 \theoremstyle{definition}
 \newtheorem{defn}[thm]{Definition}
 \theoremstyle{remark}
 
 \newtheorem*{ex}{Example}
 \numberwithin{equation}{section}

\begin{document}


\newcounter{alphabet}
\newcounter{tmp}
\newenvironment{Thm}[1][]{\refstepcounter{alphabet}%
\smallskip%
\noindent%
{\bf Theorem \Alph{alphabet}}%
\ifthenelse{\equal{#1}{}}{}{ (#1)}%
{\bf .} \itshape}{\vskip 8pt}
\newcommand{\Ref}[1]{\setcounter{tmp}{\ref{#1}}\Alph{tmp}}

\newenvironment{Lem}[1][]{\refstepcounter{alphabet}%
\smallskip%
\noindent%
{\bf Lemma \Alph{alphabet}}%
{\bf .} \itshape}{\vskip 8pt}

\newenvironment{Core}[1][]{\refstepcounter{alphabet}%
\bigskip%
\noindent%
{\bf Corollary \Alph{alphabet}}%
{\bf .} \itshape}{\vskip 8pt}

\def\be{\begin{equation}}
\def\ee{\end{equation}}

\newcommand{\blem}{\begin{lem}}
\newcommand{\elem}{\end{lem}}
\newcommand{\bthm}{\begin{thm}}
\newcommand{\ethm}{\end{thm}}
\newcommand{\bcor}{\begin{cor}}
\newcommand{\ecor}{\end{cor}}
\newcommand{\beg}{\begin{ex}}
\newcommand{\eeg}{\end{ex}}
\newcommand{\bprop}{\begin{prop}}
\newcommand{\eprop}{\end{prop}}
\newcommand{\bdefe}{\begin{defn}}
\newcommand{\edefe}{\end{defn}}
\newcommand{\bpf}{\begin{proof}}
\newcommand{\epf}{\end{proof}}

\newcommand{\IC}{{\mathbb C}}
\newcommand{\ID}{{\mathbb D}}
\newcommand{\IB}{{\mathbb B}}


%
%
%
%
%
%
%
%
%
\title[$\alpha$-Bloch spaces and Hardy spaces]
 {Landau-Bloch constants for functions in \\$\alpha$-Bloch spaces and Hardy spaces}
\author
{SH. Chen}

\address{Department of Mathematics\\
Hunan Normal University\\
Changsha, Hunan 410081\\
People's Republic of China}
\email{shlchen1982@yahoo.com.cn}

\thanks{The research was partly supported by
NSF of China (No. 11071063). 
   The research was
also partly supported by the Program for Science and Technology
Innovative Research Team in Higher Educational Institutions of Hunan
Province. Prof. X. Wang is the corresponding author. }
\author{S. Ponnusamy}
\address{Department of Mathematics\\
Indian Institute of Technology Madras\\
Chennai-600 036, India}
\email{samy@iitm.ac.in}

\author{X. Wang
}
\address{Department of Mathematics\\
Hunan Normal University\\
Changsha, Hunan 410081\\
People's Republic of China}
\email{xtwang@hunnu.edu.cn}

\subjclass{Primary: 32A10; Secondary: 32A17, 32A18, 32A35}

\keywords{Landau-Bloch constant, holomorphic function, $\alpha$-Bloch space, Hardy space.
}

\date{October 9, 2011}

\begin{abstract}
In this paper, we obtain a sharp distortion theorem for a class
of functions in $\alpha$-Bloch spaces, and as an application of it, we establish
the corresponding Landau's theorem. These results
generalize the corresponding results of Bonk, Minda and Yanagihara, and Liu,
respectively. We also prove the existence of Landau-Bloch constant
for a class of functions in Hardy spaces and the obtained result is
a generalization of the corresponding result of Chen and Gauthier.
\end{abstract}

\maketitle

\section{Introduction and main results}\label{csw-sec1}
One of the long standing open problems of determining the precise value of the
schlicht Landau-Bloch constant has attracted the attention
of many authors \cite{Al,B0, HS, L, LM,M1,M2, M-89,P}.  For holomorphic functions of several
complex variables, Landau-Bloch constant does not exist (cf. \cite{Ha, W})
unless one considers the class of functions under certain constraints, see
the works of Fitzgerald and Gong \cite{FG}, Graham and Varolin \cite{GV},
Liu \cite{LX}, and Chen and Gauthier \cite{HG1}.
The existence of the Landau-Bloch constant for  the class of holomorphic quasiregular mappings and
their related classes were investigated by Bochner \cite{Boch-46}, Hahn \cite{H},
Harris \cite{Ha}, Takahashi \cite{T} and Wu \cite{W}.

In this paper, we obtain a sharp distortion theorem (see Theorem
\ref{thm1}) for a class of holomorphic functions in $\alpha$-Bloch
spaces. As an application of Theorem \ref{thm1}, we establish the
Landau theorem for this class (see Theorem \ref{thm2}). Theorems
\ref{thm1} and \ref{thm2} generalize  the corresponding results  of
Bonk, et. al \cite{B} and Liu \cite{LX}, respectively. In Theorem
\ref{thm4}, we present the existence of the Landau-Bloch constant
for a well-known class of holomorphic functions in $H^{p}$ spaces
(cf. \cite{Duren,Z1,Z2}). Moreover Theorem \ref{thm4} is a
generalization of the corresponding result of Chen, et. al.
\cite{HG1}.



In order to state our results, we need to introduce some basic
notations. Throughout the discussion, for $b\in\mathbb{C}$, we let
$\ID(b,r)=\{z\in \mathbb{C}:\, |z-b|<r\}$. We use $\mathbb{D}$ to
denote  the open unit disk $\ID (0,1)$. Also we let
$\mathbb{C}^{n}=\{z=(z_{1},\ldots,z_{n}):\
z_{1},\ldots,z_{n}\in\mathbb{C}\}$, and for
$w=(w_{1},\ldots,w_{n})\in\mathbb{C}^{n}$,
$$\IB^n(w, r)=\left \{z\in \mathbb{C}^{n}:\; |z-w|=\sqrt{\sum^{n}_{k=1}|z_k-w_k|^2}<r\right \}
$$
and
$\IB^n=\IB^n(0, 1)$. The class of all holomorphic functions from
$\mathbb{B}^{n}$ into $\mathbb{C}^{n}$ is denoted by
${\mathcal H}(\mathbb{B}^{n})$ (cf. \cite{Z2}).
Here and in the following, we
always treat $z\in \IC^n$ as a column vector, that is, $n\times 1$ column matrix
$$z=\left(\begin{array}{cccc}
z_{1}   \\
\vdots \\
 z_{n}
\end{array}\right).
$$
Thus, for any $f=(f_{1},\ldots,f_{n})\in {\mathcal H}(\IB^n)$,
we denote by $\partial f/\partial z_{k}$ the column vector formed by $\partial
f_{1}/\partial z_{k},\ldots,\partial f_{n}/\partial z_{k}$ and by
$$f'=\left (\frac{\partial f}{\partial z_{1}},\ldots,\frac{\partial f}{\partial z_{n}}\right )
:=\left (\frac{\partial f _i}{\partial z_{j}}\right )_{n\times n}
$$
we mean the matrix formed by these column vectors, namely the $n\times n$-matrix with $(i,j)$-th entry
as $\partial f_i/\partial z_{j}$.

 For an $n\times n$ matrix $A$, the operator norm is defined by
$$|A|=\sup_{x\neq 0}\frac{|Ax|}{|x|}=\max\{|A\theta|:\
\theta\in\partial \mathbb{B}^{n}\}.
$$
Let $\mathcal{B}_{\alpha}$ denote the  \textit{$\alpha$-Bloch space}
which consists of all $f\in {\mathcal H}(\mathbb{B}^{n})$
such that (see \cite{Zhu-93}) 
$$\|f\|_{\alpha}+|f(0)|<\infty,
$$
where $\alpha >0$ and $\|f\|_{\alpha}$ denotes the \textit{$\alpha$-Bloch  semi-norm} of
$f$ defined by
$$\|f\|_{\alpha} =\sup_{z\in\mathbb{B}^{n}}(1-|z|^{2})^{\alpha}| f'(z)|.
$$
Recently, many authors investigated the properties of $\alpha$-Bloch
functions in $\mathbb{B}^{n}$, see \cite{L1,L2,ZX,Z1,Z2}.

The following lemma is easy to derive and so we omit its proof.

\begin{lem}\label{lem1}
For $x\in[0,1]$, let
$$\varphi(x)=x(1-x^{2})^{\frac{\alpha(n+1)}{2}}\sqrt{\alpha(1+n)+1}\left [\frac{\alpha(n+1)+1}
{\alpha(n+1)}\right ]^{\frac{\alpha(n+1)}{2}}
$$
and
$$a_0=\frac{1}{\sqrt{\alpha(1+n)+1}}.
$$
Then $\varphi$ is increasing in $ [0,a_0 ]$, decreasing in $ [a_0 ,1]$ and
$\varphi (a_0)=1.$
\end{lem}

As in \cite{LX}, the \textit{semi-norm $\|f\|_{0,\alpha}$} of
$f\in {\mathcal H}(\mathbb{B}^{n})$ is given by
$$\|f\|_{0,\alpha}=\sup_{z\in\mathbb{B}^{n}}
\big\{(1-|z|^{2})^{\frac{\alpha(n+1)}{2n}}|\det f'(z)|^{\frac{1}{n}}\big\}.
$$

We now state our main results and their proofs will be given in Section \ref{sec-3}.

\begin{thm}\label{thm1}
Suppose that $f\in {\mathcal H}(\mathbb{B}^{n})$ such that  $\|f\|_{0,\alpha}=1$ and $\det
f'(0)=\lambda\in(0,1]$. Then for any $z$ with
$|z|\leq \frac{a_0+m(\lambda)}{1+a_0m(\lambda)}$, we have
\be\label{eq1.07}
|\det f'(z)|\geq {\rm Re\,}\big(\det f'(z)\big)\geq
\frac{\lambda(m(\lambda)-|z|)}{m(\lambda)(1-m(\lambda)|z|)^{\alpha(n+1)+1}},
\ee
where $m(\lambda)$ is the unique real root of the equation
$\varphi(x)=\lambda$ in the interval
$[0,a_0]$ and, $\varphi$ and $a_0$ are defined as in Lemma {\rm \ref{lem1}}.

Moreover, for any $z$ with $|z|\leq\frac{a_0-m(\lambda)}{1-a_0m(\lambda)}$, we have
\be\label{eq1.08}
|\det f'(z)|\leq\frac{\lambda(m(\lambda)+|z|)}{m(\lambda)(1+m(\lambda)|z|)^{\alpha(n+1)+1}}.
\ee
The estimates  $(\ref{eq1.07})$ and $(\ref{eq1.08})$ are sharp.
\end{thm}

We remark that when $\lambda=1$ and $\alpha=1$, the inequality
(\ref{eq1.07}) coincides with \cite[Theorem 5]{LX}. In particular,
when $\alpha=n=1$, Theorem \ref{thm1} coincides with \cite[Theorem
2]{B}.

For $f\in {\mathcal H}(\mathbb{B}^{n}),$ a {\it schlicht ball} of
$f$ centered at $f(w_{1})$ is a ball with center $f(w_{1})$ such
that $f$ maps an open subset of $\mathbb{B}^{n}$ containing $w_{1}$
biholomorphically onto this ball. For a point
$w_{1}\in\mathbb{B}^{n}$, let $r(w_{1},f)$ denote the radius of the
largest schlicht ball of $f$ centered at $f(w_{1})$.

\begin{thm}\label{thm2}
Suppose that $f\in\mathcal{B}_{\alpha}$  such that
$\|f\|_{\alpha}\leq K\ \mbox{for some constant}\ K\geq1,$
$\|f\|_{0,\alpha}=1$ and $\det f'(0)=\lambda\in(0,1]$. Then
$f(\mathbb{B}^{n})$ contains a schlicht ball of radius
$$r(0,f)\geq\frac{\lambda
K^{1-n}}{m(\lambda)}\int_{0}^{m(\lambda)}
\frac{(1-t^{2})^{\alpha(n-1)}(m(\lambda)-t)}{(1-m(\lambda)t)^{\alpha(n+1)+1}}\,dt,
$$
where  $m(\lambda)$ is the same as in Theorem {\rm \ref{thm1}}.
\end{thm}

We remark that when $\lambda= \alpha=1$, Theorem \ref{thm2} coincides with \cite[Theorem 6]{LX}.

Finally, for $0<p<+\infty$, the \textit{Hardy space $H^{p}$} consists of
$f\in {\mathcal H}(\mathbb{B}^{n})$ such that
$$\|f\|_{p}=\sup_{0<r<1}\left\{\int_{\partial\mathbb{B}^{n}}|f(r\zeta)|^{p}\,d\sigma(\zeta)\right\}^{\frac{1}{p}}<+\infty,
$$
where $d\sigma$ denotes the normalized surface measure.

Now we state our final result which gives the existence of Landau-Bloch constant for a
class of holomorphic functions in Hardy spaces.

\begin{thm}\label{thm4}
Suppose that $f\in H^{p}$ satisfies $\|f\|_{p}\leq K_{0}$ for some
constant $K_0>0$, $f(0)=0$ and $|\det f'(0)|=\lambda_{0}>0$, where
$K_{0}\geq\lambda_{0}$. Then $f$ is univalent in
$\mathbb{B}^{n}(0,\rho_{1}(r_{0}))$  with
$$\rho_{1}(r_{0})=\max_{0<r<1}\rho_{1}(r),
$$
where 
$$\rho_{1}(r)=\frac{\lambda_{0}}{mK_{0}^{n}}\big[r^{n+1}(1-r^{2})^{\frac{n^{2}}{p}}\big],~r_{0}=\sqrt{\frac{p(n+1)}{p(n+1)+2n^{2}}}
$$ 
and
$$m=\frac{\sqrt{2}(7+\sqrt{17})}{4\sqrt{5-\sqrt{17}}}\approx 4.199556 .
$$
Moreover, the range
$f(\mathbb{B}^{n}(0,\rho_{1}))$ contains a univalent ball
$\mathbb{B}^{n}(0,R)$ with
$$R=\frac{\lambda_{0}^{2}}{2mK_{0}^{2n-1}}\sqrt{\frac{p(n+1)}{p(n+1)+2n^{2}}}\left(\frac{2n^{2}}
{p(n+1)+2n^{2}}\right )^{\frac{2n^{2}-n}{p}}.
$$
\end{thm}
We remark that Theorem \ref{thm4} is a generalization of
\cite[Theorem 2]{HG1}.
\section{Preliminaries and Lemmas}

We begin to recall some basic facts about the hyperbolic geometry in $\mathbb{D}$. Let
$$\rho(z,w)=\frac{1}{2}\log\left(\frac{1+|\frac{z-w}{1-\overline{z}w}|}
{1-|\frac{z-w}{1-\overline{z}w}|}\right)=\mbox{arctanh}\left |\frac{z-w}{1-\overline{z}w}\right|
$$
denote the hyperbolic distance between $z$ and $w$ in $\mathbb{D}$.
Throughout this article, we denote the hyperbolic disk (resp.
circle) by $\mathbb{D}_{h}(b,r)=\{z:\ \rho(b,z)<r\}$ (resp.
$\mathbb{S}_{h}(b,r)=\{z:\ \rho(b,z)=r\}$). Hyperbolic disks and
circles in $\mathbb{D}$ are actually Euclidean disks and circles,
respectively,  with possibly different centers and radii. For $b$
and $z\in\mathbb{D}$, we can easily see that
\be\label{xx}
\rho(b,z)=r ~\Longleftrightarrow ~ \left
|\frac{z-b}{1-\overline{b}z}\right|=\mbox{tanh}(r)
~\Longleftrightarrow ~
\frac{|1-\overline{b}z|^{2}}{1-|z|^{2}}=\frac{1-|b|^{2}}{1-\mbox{tanh}^{2}(r)}.
\ee

For the proof of our main results we need several lemmas.

\begin{Lem}{\rm (\cite[Lemma 4]{LX})}\label{LemA}
Let $A$ be an $n \times n$ complex matrix with $|A|>0$. Then for any
unit vector $\theta\in\partial \mathbb{B}^{n}$, the inequality
$$|A\theta|\geq\frac{|\det A|}{|A|^{n-1}}
$$
holds.
\end{Lem}

\begin{Lem}{\rm (\cite[Lemma 3]{HG1})}\label{lem2}
Let $f$ be a holomorphic function of $\mathbb{B}^{n}$ into
$\mathbb{C}^{n}$ such that $|f(z)|\leq M$ for $z\in\mathbb{B}^{n}$,
where $M$ is a positive constant. Then
$$|f'(z)|\leq\frac{M}{1-|z|^{2}} ~\mbox{ and }~
|\det f'(z)|\leq\frac{M^{n}}{(1-|z|^{2})^{\frac{n+1}{2}}}.
$$
\end{Lem}

\begin{Lem}{\rm (\cite[Lemma 4]{HG1})}\label{LemE}
Let $A$ be a holomorphic function of $\mathbb{B}^{n}(0,r)$ into the
space of all $n \times n$ complex matrices. If $A(0)=0$ and
$|A(z)|\leq M$ for $z\in\mathbb{B}^{n}(0,r),$ then
$$|A(z)|\leq \frac{M}{r}|z|.
$$
\end{Lem}

\begin{defn}\label{def2.1}
Let $\Phi$ and $\Psi$ be holomorphic in $\mathbb{D}$ with
$\Phi(0)=\Psi(0)$ and $\mathbb{D}^{\ast}$ be an open disk in
$\mathbb{D}$ with $0\in\mathbb{D}^{\ast}$. We say $\Phi$ is
subordinate to $\Psi$ in $\mathbb{D}^{\ast}$ relative to the origin,
written $\Phi\prec\Psi$, if there is a holomorphic function $\omega$
defined in $\mathbb{D}^{\ast}$ with
$\omega(\mathbb{D}^{\ast})\subset\mathbb{D}^{\ast}$, $\omega(0)=0$
and $\Psi\circ\omega=\Phi$. If $\mathbb{D}_{h}$ is any hyperbolic
disk (relative to hyperbolic geometry in $\mathbb{D}^{\ast}$) with
center $0$, then $\Phi\prec\Psi$ implies
$\Phi(\mathbb{D}_{h})\subset\Psi(\mathbb{D}_{h})$ since the function
$\omega$ must map $\mathbb{D}_{h}$ into itself (see
\cite[page 248]{B}). 
\end{defn}

\section{Proofs}\label{sec-3}
\paragraph{Proof of Theorem \ref{thm1}.}
For the proof of the inequality (\ref{eq1.07}), we denote the
inverse of the restriction of $\varphi$ to the interval $[0,a_0]$ by
$m: \,[0,1]\rightarrow [0,a_0],$ where $\varphi$ is the same as in
Lemma \ref{lem1}. It is not difficult to see that $m$ is increasing
with $m(0)=0,$ $m(1)=a_0$ and there is an unique $m(\lambda)\in
[0,a_0]$ such that $\varphi(m(\lambda))=\lambda.$

For convenience, we set $a=m(\lambda)$ and for a fixed $\zeta\in\partial\mathbb{B}^{n}$, we define a
function $T$ on $\ID$ by
$$T(u)=(1-au)^{\alpha(1+n)}\det f'(\zeta u),\;\, u\in\mathbb{D}.
$$
 Then $T$ is a holomorphic function in $\mathbb{D}$ and
$T(0)=\lambda.$ Now, for each
$u\in\mathbb{S}_{h}\big(a,\mbox{arctanh}\,a_0 \big),$ we have
\begin{eqnarray*}
|T(u)|&=&\frac{|1-au|^{\alpha(n+1)}}{(1-|u|^{2})^{\frac{\alpha(n+1)}{2}}}(1-|u|^{2})^{\frac{\alpha(n+1)}{2}}|\det
f'(\zeta u)|\\
&\leq&\frac{|1-au|^{\alpha(n+1)}}{(1-|u|^{2})^{\frac{\alpha(n+1)}{2}}}\\
&=&\left (\frac{1-a^{2}}{1-|\frac{a-u}{1-au}|^{2}}\right)^{\frac{\alpha(n+1)}{2}}\\
&=&(1-a^{2})^{\frac{\alpha(n+1)}{2}}\left [\frac{\alpha(n+1)+1}{\alpha(n+1)}\right]^{\frac{\alpha(n+1)}{2}},
\end{eqnarray*}
which implies that $T$ maps $\mathbb{D}_{h}\big(a,\mbox{arctanh}\, a_0 \big)$
into
$$\left\{z:\,|z|<(1-a^{2})^{\frac{\alpha(n+1)}{2}}
\left [\frac{\alpha(n+1)+1}{\alpha(n+1)}\right ]^{\frac{\alpha(n+1)}{2}}\right \}.
$$
Next, we define \
$$G_{\lambda}(u)=\frac{\lambda(a-u)}{a(1-au)}.
$$
By (\ref{xx}), we observe that 
\be\label{xx1}
\rho(a,u)=\mbox{arctanh}\,a_{0} ~\Longleftrightarrow ~
\frac{\lambda}{a}\left
|\frac{u-a}{1-\overline{a}u}\right|=\frac{a_{0}\lambda}{a}.
\ee 
On the other hand, $\varphi(a)=\lambda$ implies 
\be\label{xx2}
(1-a^{2})^{\frac{\alpha(n+1)}{2}}\left [\frac{\alpha(n+1)+1}
{\alpha(n+1)}\right]^{\frac{\alpha(n+1)}{2}}=\frac{\lambda}{a\sqrt{\alpha(n+1)+1}}.
\ee 
These observations together with (\ref{xx1}) and (\ref{xx2})
show that $G_{\lambda}$ is a M\"obius transformation which is
univalent and maps $\mathbb{D}_{h}\big(a,\mbox{arctanh}\, a_0\big)$
onto $\ID(0,r_{0})$, where
$$r_0=\frac{\lambda  a_0}{a}=\frac{\lambda}{a\sqrt{\alpha(n+1)+1}}.
$$
Also, $T(0)=\lambda=G_{\lambda}(0)$ and we obtain that
$T\prec G_{\lambda}$ on $\mathbb{D}_{h}(a,\mbox{arctanh}\, a_0)$,  where
$\prec$ denotes the subordination (see Definition \ref{def2.1}). 

Now, we fix $u\in\Big (0,\frac{a_0+a}{1 +a_0a}\Big]$ and let
$\delta_{u}=\{z:\, |z|=u\}. $ Obviously,
$$\delta_{u}\subset\mathbb{D}_{h}\big(a,\mbox{arctanh}\,a_0\big).
$$
Since $T\prec G_{\lambda}$ on
$\mathbb{D}_{h}\big(a,\mbox{arctanh}\,a_0\big)$, $T$ maps the circle
$\delta_{u}$ into the closed disk bounded by the circle
$G_{\lambda}(\delta_{u})$. We find that $G_{\lambda}(\delta_{u})$ is
a hyperbolic circle  with hyperbolic center $G_{\lambda}(0)=\lambda$
and symmetric about the real axis $\mathbb{R}$. Also, it is easy to
see that $G_{\lambda}$ is decreasing on
$\mathbb{D}_{h}\big(a,\mbox{arctanh}\,a_0\big)\cap\mathbb{R}$.
Hence, $G_{\lambda}(u)$ satisfies
$$G_{\lambda}(u)=\min\big\{\mbox{Re}\big (G_{\lambda}(z)\big ):\, z\in \delta_{u} \big\},
$$
which in turn implies that
$$\mbox{Re}\,T(u) \geq\min\{|G_{\lambda}(z)|:\, z\in\delta_{u}\}=G_{\lambda}(u)
$$
whence
$$\big|\det f'(\zeta u)\big|\geq\mbox{Re\,}\big (\det f'(\zeta u)\big )\geq
\frac{\lambda(a-|u|)}{a(1-a|u|)^{\alpha(n+1)+1}}.
$$
For each $z$ with $|z|\leq\frac{a_0+a }{ 1+a_0a}$, by taking $u=|z|$
and $\zeta=z/|z|$, we have
$$|\det f'(z)|\geq \mbox{Re\,}\big (\det f'(z)\big ) \geq
\frac{\lambda(a-|z|)}{a(1-a|z|)^{\alpha(n+1)+1}}.
$$ 
The proof of (\ref{eq1.07}) is finished.

Now we shall prove (\ref{eq1.08}). For each $y\in\Big(0,\frac{a_0-a}{ 1-a_0a}\Big), $ we set
$$\delta_{-y}=\{z:\, |z|=y\}.
$$
Then
$$\delta_{-y}\subset\mathbb{D}_{h}\big(a,\mbox{arctanh}\,a_0\big).
$$
Since $T\prec G_{\lambda}$ on $\mathbb{D}_{h}\big(a,\mbox{arctanh}\,
a_0\big)$, we see that $T$ maps the circle $\delta_{-y}$ into the
closed disk bounded by the circle $G_{\lambda}(\delta_{-y})$. It is
not difficult to see that $G_{\lambda}(\delta_{-y})$ is a
hyperbolic circle with hyperbolic center $G_{\lambda}(0)=\lambda$
and symmetric about the real axis $\mathbb{R}$. We see that
$G_{\lambda}$ is decreasing on
$\mathbb{D}_{h}\big(a,\mbox{arctanh}\,a_0\big)\cap\mathbb{R}$. Then
$$G_{\lambda}(-y)=\max\{|G_{\lambda}(u)|:\, u\in\delta_{-y}\},
$$
which yields 
$$\left|(1+ay)^{\alpha(1+n)}\det f'(-y\zeta)\right|=|T(-y)|\leq
G_{\lambda}(-y)=\frac{\lambda(a+y)}{a(1+ay)}.
$$ 
This gives
$$\left|\det
f'(-y\zeta)\right|\leq\frac{\lambda(a+y)}{a(1+ay)^{\alpha(n+1)+1}}.
$$
For each $z$ with
$|z|\leq\frac{a_0-a}{ 1-a_0a}$, by taking $y=|z|$ and
$\zeta=-z/|z|$, we conclude that
$$|\det f'(z)|\leq\frac{\lambda(a+|z|)}{a(1+a|z|)^{\alpha(n+1)+1}}.
$$ 
The proof of (\ref{eq1.08}) is finished.

Moreover, the estimates  (\ref{eq1.07}) and (\ref{eq1.08}) are
sharp. The extremal function is
$$f(z)=\left(\begin{array}{cccc}
\displaystyle
\int_{0}^{z_{1}}\frac{\lambda(a-\xi)}{a(1-a\xi)^{\alpha(n+1)+1}}\,d\xi   \\
z_2\\ \vdots \\
 z_{n}
\end{array}\right).
$$
The proof of this theorem is complete. \hfill $\Box$

\vspace{6pt}

\paragraph{Proof of Theorem \ref{thm2}.}
Let $r(0,f)$ be the supremum of the positive numbers $r$ such that
there exists a domain $\Omega_{r}\subset\mathbb{B}^{n}$ containing
$0$ which is mapped biholomorphically onto the ball
$\mathbb{B}^{n}(0,r)$. Then, since $f$ is locally biholomorphic, it
follows from the monodromy theorem that there is a point $w_{0}$ on
the boundary of the ball $\mathbb{B}^{n}(0,r(0,f))$ such that the
arc $\gamma=f^{-1}[0,w_{0}]$ originates from the origin and tends to
$\partial\mathbb{B}^{n}$ or  to a point $z_{0}\in\mathbb{B}^{n}$
with $|z_{0}|\geq m(\lambda)$ and $\det f'(z_{0})=0$ as
$w\rightarrow w_{0}$. Let $\Gamma=[0,w_{0}]$ denote the segment from
$0$ to $w_{0}.$ Because $\|f\|_{\alpha}$ is finite by hypothesis, it
follows from the definition of $\mathcal{B}_{\alpha}$ that
$$|f'(z)|\leq\frac{\|f\|_{\alpha}}{(1-|z|^{2})^{\alpha}} ~\mbox{ for $z\in\mathbb{B}^{n}$}.
$$
Therefore, by Lemma~A 
and (\ref{eq1.07}), we obtain that
\begin{eqnarray*}
r(0,f)&=&\int_{\Gamma}\,|dw| =\int_{\gamma}\left |f'(z)\frac{dz}{|dz|}\right|\cdot |dz|\\
&\geq&\int_{\gamma}\frac{|\det f'(z)|}{|f'(z)|^{n-1}}\,|dz|\\
&\geq& \frac{\lambda
K^{1-n}}{m(\lambda)}\int_{0}^{m(\lambda)}
\frac{(1-t^{2})^{\alpha(n-1)}(m(\lambda)-t)}{(1-m(\lambda)t)^{\alpha(n+1)+1}}\,dt
\end{eqnarray*}
and the proof of this theorem is complete. \hfill $\Box$

\vspace{6pt}

\paragraph{Proof of Theorem \ref{thm4}.}
For a fixed $r\in(0,1)$, let $F(z)=r^{-1}f(rz)$ for
$z\in\mathbb{B}^{n}$. For $\rho^{\ast}\in(0,1)$, it is easy 
to see that if $F$ is univalent in $\mathbb{B}^{n}(0,\rho^{\ast})$,
then $f$ is univalent in $\mathbb{B}^{n}(0,r\rho^{\ast})$. On the
other hand, for $R^{\ast}\in(0,1)$, if
$F(\mathbb{B}^{n}(0,\rho^{\ast}))$ contains an univalent ball
$\mathbb{B}^{n}(0,R^{\ast})$, then
$f(\mathbb{B}^{n}(0,r\rho^{\ast}))$ contains an univalent ball
$\mathbb{B}^{n}(0,rR^{\ast})$.

Now, we begin to prove the univalency of $F$ in $\mathbb{B}^{n}(0,\rho_{0}(r))$.
By \cite[Theorem 4.17]{Z2}, we find that
$$|F(z)|\leq\frac{\|f\|_{p}}{r(1-r^{2})^{n/p}}\leq M_{0}(r),
$$
where $M_{0}(r)=K_{0}/(r(1-r^{2})^{n/p}).$ For each $z\in\mathbb{B}^{n}$, Lemma~B 
yields that
$$|F'(z)-F'(0)|\leq M_{0}(r)\left(1+\frac{1}{
1-|z|^{2}}\right)=\frac{M_{0}(r)(2-|z|^{2})}{1-|z|^{2}}.
$$
Let $W_{1}(r)=(2-r^{2})/[r(1-r^{2})] \ \mbox{for}\ r\in(0,1).$ It is
easy to see that
$$W_{1}(r_{1})=\min_{r\in(0,1)}\{W_{1}(r)\},
$$
where $r_{1}=\sqrt{\frac{5-\sqrt{17}}{2}}\approx0.662.$ We denote
$W_{1}(r_{1})$ by $m.$ Then, we can easily verify that
$$m:=W_1(r_1)= \frac{\sqrt{2}(7+\sqrt{17})}{4\sqrt{5-\sqrt{17}}}\approx 4.199556 .
$$
For  $|z|\leq r_{1}$, by Lemma~C, 
we have
$$|F'(z)-F'(0)|\leq mM_{0}(r)|z|.
$$
On the other hand, for any $\theta\in\partial\mathbb{B}^{n}$, we infer from Lemma~A 
that
$$|F'(0)\theta|\geq\frac{\lambda_{0}}{|F'(0)|^{n-1}}\geq\frac{\lambda_{0}}{M_{0}^{n-1}(r)}.
$$

In order to prove the univalence of $F$ in
$\mathbb{B}^{n}(0,\rho_{0}(r))$ with
$\rho_{0}(r)=\frac{\lambda_{0}}{mM_{0}^{n}(r)}<r_{1}$, we choose two
distinct points $z'$, $z''\in \mathbb{B}^{n}(0,\rho_{0}(r))$ and let
$[z',z'']$ denote the segment from $z'$ to $z''$ with the endpoints
$z'$  and $z''$.
Set 
$$dz=\left(\begin{array}{cccc}
dz_{1}   \\
\vdots \\
 dz_{n}
\end{array}\right)\;\; \mbox{and}\;\; d\overline{z}=\left(\begin{array}{cccc}
d\overline{z}_{1}   \\
\vdots \\
 d\overline{z}_{n}
\end{array}\right).
$$
Then we have
\begin{eqnarray*}
|F(z')-F(z'')|&\geq&
\left|\int_{[z',z'']}F'(0)\,dz\right|
-\left|\int_{[z',z'']}(F'(z)-F'(0))\,dz\right|\\
&>&|z'-z''|\left\{\frac{\lambda_{0}}{M_{0}^{n-1}(r)}-mM_{0}(r)\rho_{0}(r)\right\}\\
&=&0
\end{eqnarray*}
which shows that $F$ is univalent in
$\mathbb{B}^{n}(0,\rho_{0}(r))$. Then $f$ is univalent in
$\mathbb{B}^{n}(0, \rho_{1}(r))$, where $\rho_{1}(r)=r\rho_{0}(r).$
By calculations, we have
\begin{eqnarray*}
\max_{0<r<1}\rho_{1}(r)=\frac{\lambda_{0}}{mK_{0}^{n}}\max_{0<r<1}\left[r^{n+1}(1-r^{2})^{\frac{n^{2}}{p}}\right]
=\rho_{1}(r_{0} ),
\end{eqnarray*} 
where 
$$r_{0}=\sqrt{\frac{p(n+1)}{p(n+1)+2n^{2}}}.
$$

Now, for  $z$ with $|z|=\rho_{0}(r_{0})=\frac{\lambda_{0}}{mM_{0}^{n}(r_{0})}$, we have
\begin{eqnarray*}
|F(z)-F(0)|&\geq& \left|\int_{[0,\rho_{0}(r_{0})]}F'(0)\,dz\right|
-\left|\int_{[0,\rho_{0}(r_{0})]}(F'(z)-F'(0))\,dz\right|\\
&\geq&\frac{\lambda_{0}\rho_{0}(r_{0})}{M_{0}^{n-1}(r_{0})}-mM_{0}(r_{0})\int_{0}^{\rho_{0}(r_{0})}\rho
d\rho\\
&=&\rho_{0}(r_{0})\left\{\frac{\lambda_{0}}{M_{0}^{n-1}(r_{0})}-\frac{mM_{0}(r_{0})\rho_{0}(r_{0})}{2}\right\}\\
&=&\frac{\lambda^{2}_{0}}{2mM_{0}^{2n-1}(r_{0})}.
\end{eqnarray*}
Hence $F\big(\mathbb{B}^{n}(0,\rho_{0}(r_{0}))\big)$ contains a
univalent ball $\mathbb{B}^{n}(0,R_{0})$ with
$$R_{0}=\frac{\lambda^{2}_{0}}{2mM_{0}^{2n-1}(r_{0})},
$$ 
which implies $f\big(\mathbb{B}^{n}(0,\rho_{1}(r_{0}))\big)$ contains a
univalent ball $\mathbb{B}^{n}(0,R)$ with

\vspace{7pt} \hfill 
$R=rR_{0}= \frac{\lambda_{0}^{2}}{2mK_{0}^{2n-1}}\sqrt{\frac{p(n+1)}{p(n+1)+2n^{2}}}\left(\frac{2n^{2}}
{p(n+1)+2n^{2}}\right )^{\frac{2n^{2}-n}{p}}. $\hfill $\Box$


\subsection*{Acknowledgment}
The authors would like to thank the referee
for his/her careful reading of this paper and many useful
suggestions.
\end{document}